\documentclass{amsart}
\usepackage{amsmath}
\usepackage{amsfonts}
\usepackage{amssymb}
\usepackage{graphicx}
\usepackage{float, verbatim}
\usepackage{enumerate}
\usepackage{color}
\usepackage{caption}
\usepackage{subcaption}
\usepackage{palatino}

\makeatletter
\@namedef{subjclassname@2010}{%
	\textup{2010} Mathematics Subject Classification}
\makeatother

\newtheorem{thm}{Theorem}[section]

\newtheorem{ques}[thm]{Question}
\newtheorem{lma}[thm]{Lemma}
\newtheorem{cor}[thm]{Corollary}
\newtheorem{defn}[thm]{Definition}

\newtheorem{conj}[thm]{Conjecture}
\newtheorem{rem}[thm]{Remark}
\newtheorem*{thm*}{Theorem}
\newtheorem*{conj*}{Conjecture}
\newtheorem*{rem*}{Remark}

\newcommand{\ds}{\textit{ Duffin-Schaeffer conjecture }}

\begin{document}

\title{A Fourier analytic approach to inhomogeneous Diophantine approximation}

\author{Han Yu}
\address{Han Yu\\
School of Mathematics \& Statistics\\University of St Andrews\\ St Andrews\\ KY16 9SS\\ UK \\ }
\curraddr{}
\email{hy25@st-andrews.ac.uk}
\thanks{}



\date{}

\dedicatory{}

\begin{abstract}
In this paper we study inhomogeneous Diophantine approximation with rational numbers of reduced form. The central object to study is the set $W(f,\theta)$ as follows,
\begin{eqnarray*}
\left\{x\in [0,1]:\left |x-\frac{m+\theta(n)}{n}\right|<\frac{f(n)}{n}\text{ for infinitely many coprime pairs of numbers } m,n\right\},
\end{eqnarray*}
where $\{f(n)\}_{n\in\mathbb{N}}$ and $\{\theta(n)\}_{n\in\mathbb{N}}$ are sequences of real numbers in $[0,1/2]$. We will completely determine the Hausdorff dimension of $W(f,\theta)$ in terms of $f$ and $\theta$. As a by-product, we also obtain a new sufficient condition for $W(f,\theta)$ to have full Lebesgue measure and this result is closely related to the study of \ds with extra conditions.
\end{abstract}

\maketitle
\section{Introduction of the results}
We are interested in Diophantine approximation with inhomogeneous shifts. Although it may look similar, the nature of inhomogeneous Diophantine approximation is considered to be rather different from its homogeneous counterpart \cite{L}. A nice introduction to the field can be found in \cite{BU}. There are also some recent results on inhomogeneous Diophantine approximation that come from different aspects of metric number theory and dynamical systems, see \cite{LN}, \cite{RA} for more details. In \cite{SC}, Chow proved a result which is closely related to the inhomogeneous Littlewood conjecture. The conjectures \cite[Conjecture 1.6, 1.7]{SC} also give some motivation for the content of this paper.

We now introduce the following sets of well approximable numbers. In the statement, we write `i.m.' for `infinitely many'.
 
\begin{defn}
	Given any sequences \[f:\mathbb{N}\to [0,1/2] \text{ and }\theta:\mathbb{N}\to [0,1/2],\] we define the following two sets
	
	\[
	W_0(f,\theta)=\left\{x\in [0,1]:\left |x-\frac{m+\theta(n)}{n}\right|<\frac{f(n)}{n} \text{ for i.m. numbers } m,n\right\}
	\]
	and
	\[
	W(f,\theta)=\left\{x\in [0,1]:\left |x-\frac{m+\theta(n)}{n}\right|<\frac{f(n)}{n} \text{ for i.m. coprime pairs of numbers } m,n\right\}.
	\]
	
	We call the sequence $f$ an approximation function and the sequence $\theta$ an inhomogeneous shift. The sets $W_0(.,.),W(.,.)$ are called sets of well approximable numbers with respect to $f,\theta$.
\end{defn} 

When $\theta$ is constantly equal to $0$ or equivalently $\theta=\mathbf{0}$, the study of well approximable numbers  is referred to as homogeneous Diophantine approximation. In this case we have a good understanding about the size (in terms of the Lebesgue measure) of $W_0(.,.)$ and some partial information about $W(.,.)$. See \cite[Chapter 2]{BV3} for some detailed discussions. However, when $\theta$ is not the zero function, we encounter inhomogeneous Diophantine approximation. So far we do not have a complete understanding about the size of $W(.,.)$ in terms of either Lebesgue or the Hausdorff dimension.
\subsection{The Hausdorff dimension is invariant under inhomogeneous shifts}

The Hausdorff dimension of $W_0(.,\mathbf{0})$ was studied extensively by Hinokuma and Shiga in \cite{HS}. In fact there is an explicit formula for computing $\dim_{H} W_0(f,\mathbf{0})$ in terms of the approximation function $f$. We will introduce this formula later. In this paper we are interested in the Hausdorff dimension of $W(.,.)$. Our main theorem in this paper is as follows. In below we use $\phi(.)$ for the Euler totient function and $d(.)$ for the divisor function, see Section \ref{Notations} for more details.

\begin{thm}\label{ThmDimMain}
	For any approximation function $f$ and inhomogeneous shift $\theta$, we have the following equality
	\[
	\dim_H W(f,\theta)=\dim_H W_0(f,\mathbf{0}).
	\]
\end{thm}

In particular the above result implies that the Hausdorff dimension of $W(f,\theta)$ depends only on $f$ and not on $\theta$. Thus we have completely determined the Hausdorff dimension of sets of well approximable numbers with reduced fractions and arbitrary inhomogeneous shift. We shall see that Theorem \ref{ThmDimMain} is a consequence of the following result. 
\begin{thm}\label{ThmBy}
	For any approximation function $f$ and inhomogeneous shift $\theta$, we have the following result
	\[
	\limsup_{N\to\infty} \frac{\sum_{n=1}^N \phi(n)f(n)/n}{\sqrt{\sum_{n=1}^N f(n)d^3(n)\log^2 n}}=\infty\implies
	W(f,\theta) \text{ has full Lebesgue measure.}
	\]
\end{thm}

The above theorem is closely related with \cite[Conjecture 1.7]{SC}. Ideally we want to get rid of the logarithmic and divisor function in the denominator. It is actually possible to prove the following result if we use all fractions instead of only the reduced ones.
\begin{thm*}[See Theorem \ref{Thmpo2} below]
	For any approximation function $f$ and inhomogeneous shift $\theta$, we have the following result
	\[
	\limsup_{N\to\infty} \frac{\sum_{n=1}^N f(n)}{\sqrt{\sum_{n=1}^N f(n)d(n)}}=\infty\implies
	W_0(f,\theta) \text{ has full Lebesgue measure.}
	\]
\end{thm*}
 In Section \ref{Chow} we shall discuss these resuts further. We note here that it is also possible to estimate the growth of the number of approximating fractions for a Lebesgue typical point. For more precise descriptions and discussions, see Theorem \ref{ThmNm} below.

If we use the result about the maximal order of the divisor function $d(.)$ we can get the following corollary which is easier to work with.

\begin{cor}\label{Thm1}
	For any approximation function $f$ and inhomogeneous shift $\theta$, if
	\[
	\limsup_{N\to\infty} \frac{\sum_{n=1}^N \phi(n)f(n)/n}{\log^2 N\log\log N\exp(3\log 2 \log N/\log\log N)}=\infty
	\]
	then
	\[
	W(f,\theta) \text{ has full Lebesgue measure.}
	\]
	To obtain an even more convenient result, we can replace the denominator with $N^{\epsilon}$ for any $\epsilon>0$.
\end{cor}

The homogeneous version of the above result, in a slightly different form, appeared in \cite{HPV} and was improved later in \cite{BV2}. We will discuss these results later in this paper.

\subsection{Some further results about inhomogeneous Diophantine approximation}
Our method can help us deal with the Lebesgue measure of $W(.,.)$ in some cases. Our next result is related with the study of the \ds with extra conditions. This topic was studied in \cite{BV2} and \cite{HPV}. With Theorem \ref{ThmBy} and Corollary \ref{Thm1} above, we can revisit \cite[Theorem 1]{HPV} for inhomogeneous Diophantine approximation and it is interesting to see how much more we can obtain for Diophantine approximation with inhomogeneous shift. Our general result is as follows. In below $\mathcal{H}^h(.)$ denotes the Hausdorff measure with dimension function $h$, more details can be found in \cite[Section 2]{BV} and the references therein.

\begin{thm}\label{ThDim}
	For any approximation function $f$ and inhomogeneous shift $\theta$, let $h:\mathbb{R}^+\to\mathbb{R}^+$ be such that $h(x)\to 0$ as $x\to 0$ and $h(x)/x$ is monotonic. If the following condition holds
	\[
	\limsup_{N\to\infty} \frac{\sum_{n=1}^N \phi(n)h(f(n)/n)}{\log^{2.5} N\left(\max_{n\in [1,N]} h(f(n)/n)^{1/2}n\right)}=\infty,
	\]
	then
	\[
   \mathcal{H}^h \left(W(f,\theta)\right)=\mathcal{H}^h([0,1]).
	\]

\end{thm}

In Section \ref{sefurther} we will provide an example to show that the above theorem is not covered by known results in the homogeneous case. The above theorem is rather complicated to use in practice and we shall obtain the following corollary which is more convenient to work with.

\begin{cor}\label{Co1}
	For any approximation function $f$ and inhomogeneous shift $\theta$, if there exists a number $A>3$ such that
	\[
	f(n)=O\left(\frac{\log ^A n}{n}\right)
	\]
	and
	\[
	\limsup_{N\to\infty}\frac{\sum_{n=1}^{N}\frac{f(n)}{n}\phi(n)}{\log ^{\frac{A}{2}+2.5}N}=\infty,\label{Con2}\tag{1}
	\]
	then
	\[
	 W(f,\theta) \text{ has full lebesgue measure.}
	\]
	The conclusion holds true if we replace the condition $(1)$ with
	\[
	\sum_{n=2}^{\infty}\frac{f(n)}{\log ^{A/2+2.5+\epsilon} n}=\infty,\tag{2}
	\]
	for some $\epsilon>0.$
	If $\theta=\mathbf{0}$, then condition $(1)$ can be weakened slightly because of a result of Gallagher \cite{GA} to the following
	\[
	\limsup_{N\to\infty}\frac{\sum_{n=1}^{N}\frac{f(n)}{n}\phi(n)}{\log ^{\frac{A}{2}+2.5}N}>0.\tag{1'}
	\]
\end{cor}

\begin{proof}
	We set $h(x)=x$ in Theorem \ref{ThDim} and the first conclusion is easy to see. For the second conclusion, we assume condition $(2)$ and consider the iterated exponential intervals \[I_k=\left[2^{2^k},2^{2^{k+1}}\right].\] There are infinitely many $k>0$ such that
	\[
	\sum_{n\in I_k}\frac{f(n)}{\log ^{A/2+2.5+\epsilon}n}>\frac{1}{k^2},
	\]
	otherwise the following sum
	\[
	\sum_{n=2}^{\infty} \frac{f(n)}{\log ^{A/2+2.5+\epsilon}n}
	\]
	will not diverge.
	Then, we see that
	\[
	\sum_{n\in I_k}\frac{f(n)}{n}\phi(n)=\sum_{n\in I_k}\frac{f(n)}{\log^{A/2+2.5+\epsilon}n}\frac{\phi(n)}{n}\log^{A/2+2.5+\epsilon} n
	\]
	\begin{eqnarray*}
		&\geq& \min_{n\in I_k}\frac{\phi(n)}{n}\log^{A/2+2.5+\epsilon} n \sum_{n\in I_k}\frac{f(n)}{\log^{A/2+2.5+\epsilon}n}\\
		&\geq& C\frac{1}{\log\log 2^{2^{k+1}}}\log^{A/2+2.5+\epsilon}2^{2^{k}} \frac{1}{k^2}\\
		&\geq& C'\frac{1}{k^3} 2^{(A/2+2.5+\epsilon)k}\\
		&\geq& C'' 2^{(A/2+2.5+0.5\epsilon)(k+1)},\\
	\end{eqnarray*}
	where $C,C',C''$ are constants which only depend on $A,\epsilon$. The choice of constant $C$ comes from the following well-known result concerning  the Euler gamma $\gamma$:
	\[
	\liminf_{n\to\infty}\frac{\phi(n)}{n}\log\log n=e^{-\gamma}.\tag{***}\label{FACT}
	\]
	Then, for all $k>0$ we have
	\[
	\sum_{n=2}^{2^{2^{k+1}}}\frac{f(n)}{n}\phi(n)\geq \sum_{n\in I_k}\frac{f(n)}{n}\phi(n)\geq C'' 2^{(A/2+2.5+0.5\epsilon)(k+1)}.
	\]
	This implies condition $(1)$ because $2^{0.5\epsilon k}\to\infty$ as $k\to\infty$.
\end{proof}

\section{Some earlier results and discussions}\label{Se4}

Before the proofs we shall briefly introduce some known results in metric Diophantine approximation and discuss how our results can be related to them.  In Section \ref{sefurther} we will also discuss some related questions.

\subsection{Some general historical remarks}
One of the most famous result was first proved by Khintchine and generalized by Groshev, see for example \cite[Theorem 1]{BHH17}.
\begin{thm*}[KG]
	If $f$ is a non-increasing approximation function such that
	\[
	\sum_{n=1}^{\infty} f(n)=\infty,
	\]
	then for any $a\in [0,1/2]$ 
	$
	W_0(f,\mathbf{a}) \text{ has full Lebesgue measure.}
	$
\end{thm*}
For convenience, the bold letter $\mathbf{a}$ denotes the constant sequence $\theta$ such that $\theta(n)=a$ for all $n$. Later Duffin-Schaeffer \cite{Du} generalized Khintchine's result in the homogeneous case.
\begin{thm*}[DS]
	For any approximation function $f$, if
	\[
	\limsup_{N\to\infty} \frac{\sum_{n=1}^{N}\frac{f(n)\phi(n)}{n}}{\sum_{n=1}^{N}f(n)}>0,\label{Con1}\tag{3}
	\]
	then
	\[
	\sum_{n=1}^{\infty} \frac{f(n)}{n}\phi(n)=\infty\implies
	W(f,\mathbf{0}) \text{ has full Lebesgue measure.}
	\]
\end{thm*}
Theorem \ref{ThmBy} and Theorem \ref{ThDim} are two inhomogeneous versions of the above result. Duffin and Schaeffer also asked whether the condition (\ref{Con1}) can be dropped. They made the following famous conjecture.
\begin{conj*}[DS]
	For any approximation function $f$ we have the following result
	\[
	\sum_{n=1}^{\infty} \frac{f(n)}{n}\phi(n)=\infty\implies
	W(f,\mathbf{0}) \text{ has full Lebesgue measure.}
	\]
\end{conj*}

\subsection{\ds with extra conditions: known results before \cite{ALMTZ18}}
A lot of work has been done since the birth of the above conjecture. Various replacements of condition (\ref{Con1}) have been found and we think that the mathoverflow webpage \cite{Online} gives a nice and brief overview. Notably, the first result of this topic with extra divergence appeared in \cite[Corollary 1]{HPV} as follows

\begin{thm*}[HPV]
	For any approximation function $f$ that satisfies the following divergence condition with a positive $\epsilon>0$
	\[
	\sum _{n=1}^\infty \frac{\phi (n) f (n)/n}{ n^{\epsilon}}= \infty,
	\]
	the set $W(f,\mathbf{0})$ has full Lebesgue measure.
\end{thm*}

In fact the $n^\epsilon$ in the denominator can be replaced by $\exp(c \log n/\log\log n)$ with a suitable constant $c>0$. This was the content of \cite[Theorem 1]{HPV}. Note that Corollary \ref{Thm1} is an inhomogeneous version of the above result. Later, in \cite{BV2}, the above result was improved to the following.

\begin{thm*}[BHHV]
	For any approximation function $f$ that satisfies the following divergence condition
	\[
	\sum _{n=1}^\infty \frac{\phi (n) f (n)/n}{ \exp (c(\log \log n)(\log \log \log n))}= \infty,\label{BV}\tag{4}
	\]
	the set $W(f,\mathbf{0})$ has full Lebesgue measure.
\end{thm*}
Our motivation for Theorem \ref{ThDim} was to replace the above condition $(4)$ with the following
\[
\sum _{n=1}^\infty \frac{\phi (n) f (n)/n}{ \log^c n}= \infty.
\]
We are not able to achieve this without the following extra upper bound \footnote{A few months after the first public version of this paper, it was proven \cite{ALMTZ18} that  this upper bound condition can be dropped for homogeneous cases. For inhomogeneous cases, it is not known whether one can get rid of this upper bound condition.}
\[
f(n)=O(\log^c n/n).
\]
We remark that in \cite{VA}, it was shown that if $f(n)=O(1/n)$ then
\[
\sum_{n=1}^{\infty} \frac{f(n)}{n}\phi(n)=\infty\implies
W(f,\mathbf{0}) \text{ has full Lebesgue measure.}
\]
We see that $W(f,\mathbf{0})$ has full Lebesgue measure if either we have a strong extra divergence condition like $(4)$, or we have a strong upper bound condition like $f(n)=O(1/n)$ together with a weak divergence condition. Theorem \ref{ThDim} and Corollary \ref{Co1} show that we can also balance the strength of the upper bound and divergence conditions. We note here that our result holds in the inhomogeneous case as well. For convenience, we introduce the following notations.

\begin{defn}
	Given two non negative numbers $a,b$, we call the condition $C(a,b)$ to be the following two conditions on the approximation function $f$,
	\[
	f(n)=O\left(\frac{\log ^a n}{n}\right),
	\]
	\[
	\limsup_{N\to\infty}\frac{\sum_{n=1}^{N}\frac{f(n)}{n}\phi(n)}{\log ^{b}N}=\infty.
	\]
	We say that the condition $C(a,b)$ is sufficient if the set $W(f,\mathbf{0})$ has full Lebesgue measure under the condition $C(a,b)$.
\end{defn}
Thus the result in \cite{VA} that was mentioned earlier says that the condition $C(0,0)$ is sufficient. Corollary \ref{Co1} says that the condition $C(A,A/2+2.5+\epsilon)$ is sufficient for any $A>3,\epsilon>0$. However, it is easy to check that if the condition $C(a,b)$ is sufficient, then for any positive number $c$, the condition $C(a+c,b+c)$ is sufficient as well. Indeed, if we assume the condition $C(a+c,b+c)$ we can find the following new approximation function
\[
f'(n)=\frac{f(n)}{\log^c n}.
\]
It is then easy to check that $f'$ satisfies the condition $C(a,b)$ and because it is sufficient we see that $
W(f',\mathbf{0})
$ has full Lebesgue measure. It is clear that $f'(n)\leq f(n)$ for all sufficiently large integer $n$ and therefore $W(f,\mathbf{0})$ has full Lebesgue measure. Similarly if the condition $C(a,b)$ is sufficient then so is the condition $C(a',b)$ and $C(a,b')$ for all $0\leq a'\leq b$ and $b'\geq b$. We plot the following graph to indicate the known and new sufficient conditions $C(a,b)$ represented as points in the Euclidean plane.
\begin{figure}[h]
	\includegraphics[width=5cm]{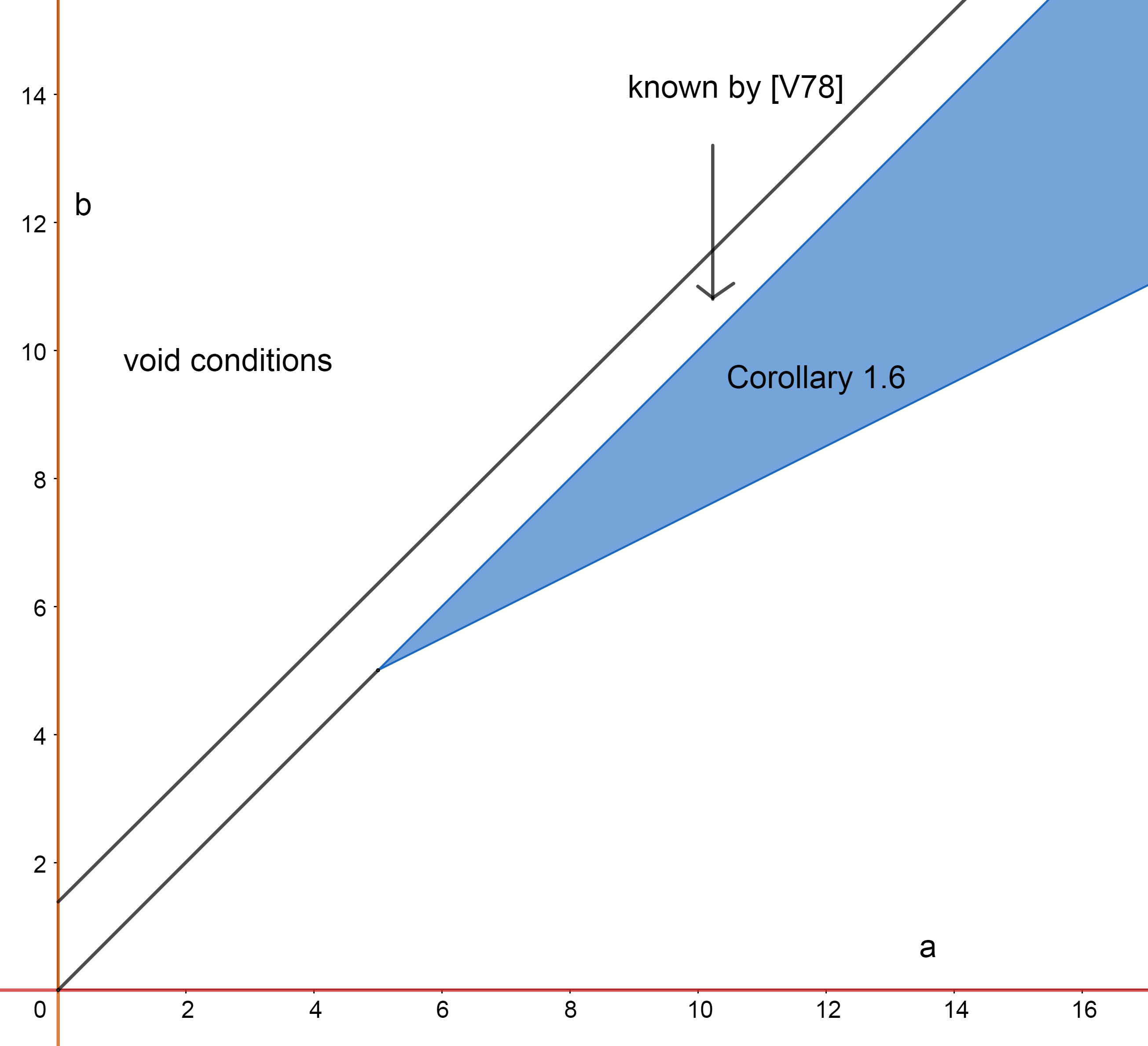}
	\caption{This picture shows sufficient conditions for $W(f,\mathbf{0})$ to have full Lebesgue measure. We remark that in the blue area, $W(f,\theta)$ also has full Lebesgue measure.}
	\label{arc2}
\end{figure}
We only considered the condition $C(a,b)$ without worrying about whether it can be satisfied at all. In fact, it is easy to check that the condition $C(a,b)$ can be satisfied if and only if $b<a+1$. This is the reason for requiring $A>3$ in Corollary \ref{Co1}.
\subsection{The Hausdorff dimension results}\label{HDR}
The Hausdorff dimensions of sets $W_0(.,.)$,$W(.,.)$ are much better known. For example we have the following result from \cite{HPV}.

\begin{thm*}[HPV]
	For any approximation function $f$ and any positive real number $s\in (0,1]$ we have the following result
	\[
	\sum_{n=1}^{\infty} \left(\frac{f(n)}{n}\right)^s\phi(n)=\infty\implies
	\dim_{H}W(f,\mathbf{0})\geq s.
	\]
\end{thm*}
By Theorem \ref{ThmDimMain}, we see that under the same condition as above, $\dim_{H}W(f,\theta)\geq s$ for all inhomogeneous shifts $\theta$. We now introduce a result for computing $\dim_H W_0(f,\mathbf{0})$ from \cite{HS}.

\begin{thm*}[HS]
	For any approximation function $f$ and real number $\alpha\in [1,\infty)$ we set
	\[
	C_\alpha(N)=\text{Cardinality of the set }\left\{n\leq N:f(n)/n\geq \frac{1}{n^{\alpha}}\right\},
	\]
	and
	\[
	\delta(\alpha)=\sup\left\{\delta:\limsup_{N\to\infty} \frac{C_{\alpha}(N)}{N^{\delta}}>0   \right\}.
	\]
	Then
	\[
	\dim_{H}W_0(f,\mathbf{0})=\min\{1,\sup_{\alpha\geq 1}\kappa(\alpha)\},
	\]
	where $\kappa(\alpha)$ is the following number
	\[
	\kappa(\alpha)=\begin{cases} 
	\frac{1+\delta(\alpha)}{\alpha}& \lim_{N\to\infty} C_{\alpha}(N)=\infty \\
	0 & \text{otherwise}
	\end{cases}
	\]
\end{thm*}

We shall see that the above theorem of Hinokuma and Shiga plays an important role in proving Theorem \ref{ThmDimMain}. As a special case, we assume now that $f$ is a non-increasing sequence and define the following lower order of $f$
\[
\lambda(f)=\liminf_{n\to\infty} \frac{-\log f(n)}{\log n}.
\]
Then for any $\alpha>\lambda(f)+1$, there exists infinitely many $n$ such that
\[
f(n)\geq \frac{1}{n^{\alpha-1}}.
\]
It follows that $\kappa(\alpha)=(1+\delta(\alpha))/\alpha$. Now, because $f$ is non increasing, we see that for any $n$ such that $f(n)\geq 1/n^{\alpha-1}$ we have
\[
f(\lfloor n/2 \rfloor)\geq\dots\geq f(n)\geq \frac{1}{n^{\alpha-1}}.
\]
As there are infinitely many such $n$, we see that for any $\alpha'>\alpha$ there are infinitely many $N$ such that
\[
C_{\alpha'}(N)\geq \frac{N}{2}.
\]
Because numbers $\alpha',\alpha$ such that $\alpha'>\alpha>\lambda(f)+1$ can be chosen arbitrarily, we see that
\[
\sup_{\alpha\geq 1}\kappa(\alpha)\geq \frac{2}{\lambda(f)+1}.
\]
On the other hand, for any $\alpha<\lambda(f)+1$, there are at most finitely many $n$ with 
\[
f(n)\geq\frac{1}{n^{\alpha-1}},
\]
and therefore we see that
\[
\kappa(\alpha)=0.
\]
This means that
\[
\sup_{\alpha\geq 1}\kappa(\alpha)=\frac{2}{\lambda(f)+1},
\]
and therefore by Theorem \ref{ThmDimMain} we see that for any inhomogeneous shift $\theta$
\[
\dim_H W(f,\theta)=\min\left\{1,\frac{2}{\lambda(f)+1}\right\}.
\]
This particular result was obtained with $W(f,\theta)$ replaced by $W_0(f,\theta)$ in \cite{L}. We note here that in \cite{L} general higher dimensional results were obtained as well.

\section{Notation}\label{Notations}

\begin{flushleft}

\begin{itemize}	
\item[1.] In this paper we always use $f$ to denote approximation functions and $\theta$ to denote inhomogeneous shifts. Unless explicitly mentioned otherwise, we assume that $f$ and $\theta$ take values in $[0,1/2]$.
\item[2.] For any number $a\in\mathbb{R}$ we use $\mathbf{a}$ to denote the constant sequence whose terms are equal to $a$.
\item[3.] We use $\dim_H$ for the Hausdorff dimension and $\mathcal{H}^h$ for the $h$-Hausdorff measure with dimension function $h$. We will not directly deal with definitions of the Hausdorff measure/dimension. For more details see \cite[Chapter 3]{Fa} and \cite[Chapter 4]{Ma}.
\item[4.] In this paper we use $\log n$ for the natural logarithm function. There is a small issue we could encounter. For an expression like $\log \log n$, we know that it is not defined at $n=1.$ Since all the results and arguments we have here deal with only the situation for $n\to\infty$, there is no problem if we simply re-define $\log 0=\log 1=2.$  
\item[5.] We shall use the following arithmetic functions:
\begin{itemize}
	\item[5.1]: The Euler function: For $n\in\mathbb{N}$,
	
	$\phi(n)=\text{ number of natural numbers smaller than and are coprime to $n$}$.
	\item[5.2]: The greatest common divisor function: For $a,b\in\mathbb{N}$
	
	$(a,b)=\text{ the greatest common divisor of $a,b$}$.
	\item[5.3]: The divisor function: For $n\in\mathbb{N},\alpha\in\mathbb{R}$,
	
	$d(n)=\text{ the number of divisors of $n$}$.
	
	\item[4]: The M\"{o}bius function:
	
		For $n\in\mathbb{N}$, the M\"{o}bius function is defined as follows,
		
		\[ 
		\mu(n)=
		\begin{cases} 
		 1 & \text{$n$ is squarefree with even number of prime factors}\\
		 -1 & \text{$n$ is squarefree with odd number of prime factors}\\
		 0 & \text{$n$ is not squarefree}
		\end{cases}
		\]
		\item[5]: The Ramanujan sum:
		
		For $n,k\in\mathbb{N}$: $c_n(k)=\sum_{1\leq a\leq n, (a,n)=1} e^{2\pi i \frac{ak}{n}}=\mu\left(\frac{n}{(k,n)}\right)\frac{\phi(n)}{\phi\left(\frac{n}{(k,n)}\right)}$
		
\end{itemize}

\item[6.] We use $P$ for general probability measure on a probability space $\Omega$ and $\lambda$ for Lebesgue measure on $[0,1]$. 

\item[7.] For a sequence of sets $A_n\subset X$: $\limsup_{n\to\infty} A_n=\{x\in X : x\in A_n \text{ for infinitely many $n\in\mathbb{N}$}\}$.

\end{itemize}
\end{flushleft}
\section{Results that will be used without proof}
The central idea we shall use in this paper is a Fourier analytic method introduced by LeVeque in \cite{LeV}. To start with, given a function $f:[0,1]\to\mathbb{R}$ which is in $L^2$ and thus in $L^1$ as well, the Fourier series of $f$ is given by
\[
\forall k\in\mathbb{N}, \hat{f}(k)=\int_0^1 e^{2\pi i kx}f(x)dx.
\]

We will need the following facts:
\[
\hat{f}(0)=\|f\|_{L^1} \text{ whenever $f$ is non negative,}
\]
\[
\widehat{fg}(0)=\sum_{k=-\infty}^{\infty} \hat{f}(k)\hat{g}(-k) \text{ whenever $f,g$ are $L^2$ functions}.
\]
The above results can be found in most text books on harmonic analysis for example in \cite[chapter 1, section 5.5]{Ka}.

We specify the version of Borel-Cantelli lemma (see \cite[lemma 2.2]{BV3}) which will be used later.
\begin{thm}\label{T1}
Let $A_n$ be a sequence of events in a probability space $(\Omega,P)$ such that
\[
\sum_{n=1}^{\infty} P(A_n)=\infty,
\]
then
\[
P(\limsup_{n\to\infty} A_n)\geq \limsup_{m\to\infty} \frac{(\sum_{n=1}^m P(A_n))^2}{\sum_{n_1,n_2=1}^{m}P(A_{n_1}\cap A_{n_2})}.
\]
\end{thm}
\begin{rem}\label{Re2}
For homogeneous metric Diophantine approximation, to conclude the full measure result we only need to show
\[
\limsup_{n\to\infty} \frac{(\sum_{n=1}^m P(A_n))^2}{\sum_{n_1,n_2=1}^{m}P(A_{n_1}\cap A_{n_2})}>0.
\]
This follows from a result of Gallagher \cite{GA}.
\end{rem}

In order to prove the general result Theorem \ref{ThDim}, we will also use the following version of the mass transference principle in \cite{BV}.
\begin{thm}[BV]\label{ThMT}
	Let $\{B_i\}_{i\in\mathbb{N}}$ be a countable collection of balls in $\mathbb{R}$ with $r(B_i)\to 0$ as $i\to \infty$. Let $h$ be a dimension function such that $h(x)/x$ is monotonic and suppose that for any ball $B$ in $\mathbb{R}$
	\[
	\lambda(B\cap\limsup_{i\to\infty} B^h_i)=\lambda(B).
	\]
	Then, for any ball $B$ in $\mathbb{R}$
	\[
	\mathcal{H}^h(B\cap\limsup_{i\to\infty} B_i)=\mathcal{H}^h(B).
	\]
	
	Here $B^h$ denotes the dilated ball. To be precise, let $B$ be a ball centred at $x\in\mathbb{R}$ with radius $r>0$ then $B^h$ is the ball centred at $x$ with radius $h(r)$.
\end{thm}
\section{Some asymptotic results on arithmetic functions}\label{Se2}
In what follows, we will use some results about the Ramanujan sum. The following result is standard and can be found in \cite{Ha} chapter 16. For integers $n,k$ we have
\[
c_n(k)=\sum_{(a,n)=1}e^{2\pi i \frac{a}{n}k}=\mu\left(\frac{n}{(n,k)}\right)\frac{\phi(n)}{\phi\left(\frac{n}{(n,k)}\right)}.
\]
We will now state and prove some technical lemmas that will be used later.
\begin{lma}\label{L1}
	There is a constant $C>0$ such that for any integers $k,m >0$ 
	\[
	\frac{1}{d(k)\log m} \sum_{n=1}^m\frac{|c_n(k)|}{\phi(n)}<C.
	\]
	Here $d(k)$ is the divisor function, that is, the number of divisors of the integer $k$.
\end{lma}
\begin{proof}
	By properties of the Ramanujan sum and the Euler totient function
	\begin{eqnarray*}
	\sum_{n=1}^m\frac{|c_n(k)|}{\phi(n)}&=& \sum_{n=1}^{m}\frac{\left|\mu\left(\frac{n}{(n,k)}\right)\right|}{\phi\left(\frac{n}{(n,k)}\right)}\\
	&=& \sum_{n=1}^{m}\frac{\left|\mu\left(\frac{n}{(n,k)}\right)\right|}{\frac{n}{(n,k)}\prod_{r|\frac{n}{(n,k)}, r\text{ prime}}(1-\frac{1}{r})}\\
	&=& \sum_{n=1}^{m}\frac{\left|\mu\left(\frac{n}{(n,k)}\right)\right|}{\prod_{r|\frac{n}{(n,k)}, r\text{ prime}}(r-1)}\\
	&=& \sum_{l=1, l\text{ squarefree}}^{m}\prod_{r|l, r\text{ prime}}\frac{1}{r-1}\left|\left\{n\in [1,m]|l=\frac{n}{(n,k)}\right\}\right|.
	\end{eqnarray*}
	The cardinality of the set can be bounded by
	\[
	\left|\left\{n\in [1,m]|l=\frac{n}{(n,k)}\right\}\right|\leq d(k),
	\]
	because $(n,k)$ must be a divisor of $k$, and for every such divisor $s|k$, the value of $n$ (if exists) can be uniquely determined by $sl$. Then we see that
	\begin{eqnarray*}
		\sum_{n=1}^m\frac{|c_n(k)|}{\phi(n)}&\leq&	d(k)\sum_{l=1, l\text{ squarefree}}^{m}\prod_{r|l, r\text{ prime}}\frac{1}{r-1}\\
		&\leq& d(k)\prod_{r\leq m, r\text{ prime}}\left(1+\frac{1}{r-1}\right).
	\end{eqnarray*}
Then this lemma follows by Mertens' third theorem.
\begin{thm*}[Mertens] We have
	\[\lim_{m\to\infty} \frac{1}{\log m}\prod_{r\leq m, r\text{ prime}}\left(1+\frac{1}{r-1}\right)=e^{\gamma},\]
\end{thm*}
where $\gamma$ is the Euler gamma $\gamma\approx 0.5772156$.
\end{proof}

\begin{lma}\label{L2}
	There exists a constant $C>0$ such that for all integers $n>1$,
	\[
	\sum_{k=1}^{n}\frac{d^2(k)}{k}<C\log ^3 n.
	\]
\end{lma}
\begin{proof}
	First, observe that
	\[
	d^2(k)=\sum_{l|k}d(l^2).
	\]
	Indeed for any integer with prime factorization $k=p_1^{a_1}\dots p_k^{a_n}$ we have that
	\[
	d(k)=\prod_{i=1}^{i=n}(a_i+1).
	\]
	It follows that:
	\begin{eqnarray*}
	\sum_{l|k}d(l^2)&=&\sum_{0\leq b_i\leq a_i,i\in\{1,2\dots n\}} \prod_{i=1}^{i=n}(2b_i+1)\\
	&=& \prod_{i=1}^{n}\left(\sum_{b_i=0}^{b_i=a_i} (2b_i+1)\right)\\
	&=& \prod_{i=1}^{n}(a_i+1)^2=d^2 (k).
	\end{eqnarray*}
	Then we have the following estimate
	\begin{eqnarray*}
	\sum_{k=1}^{n}\frac{d^2(k)}{k}&=&\sum_{l=1}^{n}d(l^2)\sum_{k:l|k}^{k\leq n}\frac{1}{k}\\
	&\leq&\sum_{l=1}^{n}\frac{d(l^2)}{l} (\log n+1)\\
	&\leq&(\log n+1)\sum_{m=1}^{n^2}\sum_{l:m|l^2}^{l\leq n}\frac{1}{l}\\
	&\leq&(\log n+1)\sum_{m=1}^{n^2}\sum_{l:m|l,l\in [1,n^2]}\frac{1}{l}\\
	&\leq&(\log n+1)\sum_{m=1}^{n^2}\frac{1}{m}\left(\log n^2+1\right)\\
	&\leq&(\log n+1)^2(2\log n+1)\leq C\log^3 n,
	\end{eqnarray*}
for a suitable constant $C>0$.
\end{proof}

\begin{lma}\label{L3}
	There exists a constant $C>0$ such that for any positive integer $m$,
	\[
	\sum_{1\leq n\leq m} d(n)d(m)(n,m)\leq C d^3(m) m\log m .
	\]
\end{lma}
\begin{proof}
	\begin{eqnarray*}
	\sum_{1\leq n\leq m} d(n)d(m)(n,m)&=& d(m)\sum_{r|m} \sum_{n\leq m, (n,m)=r} r d(n)\\
	&=& d(m)\sum_{r|m} r \sum_{a\leq n/r, (a,m/r)=1} d(ar)\\
	&\leq& d(m)\sum_{r|m} r\sum_{a\leq n/r, (a,m/r)=1}d(a)d(r)\\
	&\leq& C d(m) \sum_{r|m} rd(r) \frac{m}{r}\log \frac{m}{r}\\
	&\leq& C m d(m)\log m \sum_{r|m} d(r)\\
	&\leq& C m d(m)\log m \sum_{r|m} d(r^2)\\
	&=& C m d^3(m) \log m.
	\end{eqnarray*}
Here we used Dirichlet theorem for the divisor summatory function (for the constant $C$) and the first part of the proof of lemma \ref{L2}.
\end{proof}

\section{Fourier series and Diophantine approximation}\label{Se1}

From the Borel-Cantelli lemma (Theorem \ref{T1}), we see that it is important to show some properties of the measure of intersections. Now we are going to set up the Fourier analysis method.

Let $f,\theta$ be as mentioned above, we denote \[\epsilon_n=\frac{f(n)}{n}.\] Then we define the function 
\[
g_n(x):[0,1]\to \{0,1\}
\]
via the following relation
\[
g_n(x)=1\iff \left|x-\frac{m+\theta(n)}{n}\right|< \frac{f(n)}{n}, \text{ for an integer }m \text{ with }(m,n)=1. 
\]
It is clear that $g_n(x)$ is just the characteristic function on the set $A_n$, namely,
\[
A_n=\left\{x\in[0,1]| \exists 1\leq m\leq n, (m,n)=1, \left|x-\frac{m+\theta(n)}{n}\right|<\frac{f(n)}{n}\right\}.
\]
In our case $f(n)\leq\frac{1}{2}$ and therefore $A_n$ is a union of $\phi(n)$ many equal length disjoint intervals. The Lebesgue measure of $A_n$ is
\[
\|g_n\|_{L^1}=2\epsilon_n \phi(n).
\]
Now we see that $\lambda(A_n\cap A_m)=\|g_n g_m\|_{L^1}$. We need only to compute the case $n\neq m$ since otherwise the case is trivial.  By using Fourier series we can write the $L^1$-norm as
\[
\|g_n g_m\|_{L^1}=\sum_{k=-\infty}^{\infty} \hat{g_n}(k)\hat{g_m}(-k).
\]
The above equality holds whenever the series is absolutely convergent. This happens whenever $g_n,g_m$ are both $L^1$ functions. This is the case in our situation. Now we need to evaluate the Fourier series of $g_n$, it is easy to see that $g_n$ is just the characteristic function of 
\[
[-\epsilon_n,\epsilon_n]=\left[-\frac{f(n)}{n},\frac{f(n)}{n}\right]
\]
convolved with a sum of Dirac deltas
\[
\sum_{(a,n)=1}\delta\left(\frac{a+\theta(n)}{n}\right).
\]
We can also compute the Fourier series directly for $k\neq 0$
\begin{eqnarray*}
\int_0^1 e^{2\pi i kx} g_n(x)dx&=&\sum_{(a,n)=1} \int_{\frac{a+\theta(n)}{n}-\epsilon_n}^{\frac{a+\theta(n)}{n}+\epsilon_n}e^{2\pi i kx}dx\\
&=& \sum_{(a,n)=1}\frac{1}{\pi k} \sin(2\pi\epsilon_n k)e^{2\pi i\frac{a+\theta(n)}{n}k}\\
&=& \frac{1}{\pi k} \sin(2\pi\epsilon_n k)c_n(k)e^{2\pi i \frac{\theta(n)}{n} k},
\end{eqnarray*}
where $c_n(k)=\sum_{(a,n)=1}e^{2\pi i \frac{a}{n}k}$is the Ramanujan sum. For $k=0$, $\hat{g}_n(0)$ is simply $\|g_n\|_{L^1}$. Hence we can express $\lambda(A_n\cap A_m)$ with the following series
\begin{eqnarray*}
\lambda(A_n\cap A_m)&=&4\epsilon_n\epsilon_m\phi(n)\phi(m)\\
&+&\frac{2}{\pi^2}\sum_{k=1}^{\infty}\frac{\sin(2\pi\epsilon_n k)c_n(k)\sin(2\pi\epsilon_m k)c_m(k)\cos\left(2\pi  \left(\frac{\theta(n)}{n}-\frac{\theta(m)}{m}\right) k\right)}{k^2},
\end{eqnarray*}
where we have used the fact that the values of the Ramanujan sum are real numbers and for all pairs of integers $n,k$
\[
c_n(k)=c_n(-k).
\]
We see that inhomogeneous shifts $\theta$ create just an extra $\cos(.)$ term whose modulus is bounded by $1$.

\section{proof of Theorem \ref{ThmBy} and Corollary \ref{Thm1} }\label{sepp}
It follows from the arguments in previous section that
\[
\lambda(A_n\cap A_m)\leq 4\epsilon_n\epsilon_m\phi(n)\phi(m)+\frac{2}{\pi^2}\sum_{k=1}^{\infty}\frac{|\sin(2\pi\epsilon_n k)c_n(k)\sin(2\pi\epsilon_m k)c_m(k)|}{k^2}.
\]
The basic strategy is to split the sum over $k$ up to a number $M$ which will be determined later
\[
\sum_{k=1}^{\infty}=\sum_{k=1}^{M}+\sum_{k=M+1}^{\infty}.
\]

For the first part, we use the fact that $|\sin(x)|\leq \min\{|x|,1\}$ for all $x\in\mathbb{R}$,
\begin{eqnarray*}
& &\sum_{k=1}^{M}\frac{|\sin(2\pi\epsilon_n k)c_n(k)\sin(2\pi\epsilon_m k)c_m(k)|}{k^2}\\
&\leq&\sum_{k=1}^{M}\frac{\min\{2\pi\epsilon_nk,1\}\min\{2\pi\epsilon_mk,1\}|c_n(k)c_m(k)|}{k^2}\\
&\leq& 2\pi\sum_{k=1}^{M} \frac{1}{k}\min\{\epsilon_n,\epsilon_m\} |c_n(k)||c_m(k)|.
\end{eqnarray*}

Recalling the formula for the Ramanujan sum
\[
c_n(k)=\mu\left(\frac{n}{(n,k)}\right)\frac{\phi(n)}{\phi\left(\frac{n}{(n,k)}\right)},
\]
we see that there exists an absolute constant $C>0$ satisfying the following inequality
\begin{eqnarray*}
& &\sum_{k=1}^{M}\frac{|\sin(2\pi\epsilon_n k)c_n(k)\sin(2\pi\epsilon_m k)c_m(k)|}{k^2}\\
&\leq& 2\pi\sum_{k=1}^{M} \frac{1}{k}\min\{\epsilon_n,\epsilon_m\} |c_n(k)||c_m(k)|\\
&=& 2\pi\sum_{k=1}^{M} \frac{1}{k}\min\{\epsilon_n,\epsilon_m\} \left|\mu\left(\frac{n}{(n,k)}\right)\frac{\phi(n)}{\phi\left(\frac{n}{(n,k)}\right)}\mu\left(\frac{m}{(m,k)}\right)\frac{\phi(m)}{\phi\left(\frac{m}{(m,k)}\right)}\right|\\
&\leq& 2\pi\sum_{k=1}^{M}\frac{1}{k} \min\{\epsilon_n,\epsilon_m\} (n,k)(m,k)\left|\mu\left(\frac{n}{(n,k)}\right)\mu\left(\frac{m}{(m,k)}\right)\right|\\
&\leq& 2\pi\sum_{k=1}^M \frac{1}{k}\min\{\epsilon_n,\epsilon_m\} (n,k)(m,k)\\
&\leq& C \log M d(n)d(m)(n,m)\min\{\epsilon_n,\epsilon_m\}.
\end{eqnarray*}
Here we used the fact that $\phi(n)=n\prod_{r|n,r \text{ prime}}\frac{r-1}{r}$. For the last step we see that for any divisor $s_n$ of $n$ and $r_m$ of $m$, we can sum those $k$ such that
\[
(n,k)=s_n, (m,k)=r_m.
\]
Such $k$ must be a multiple of $[s_n,r_m]$ and therefore we obtain the following result
\[
\sum_{k:(n,k)=s_n, (m,k)=r_m}\frac{1}{k} (n,k)(m,k)=\sum_{l: l\leq M/[s_n,r_m]}\frac{s_nr_m}{l[s_n,r_m]}\leq C\log M (s_n,r_m).
\]
The previous estimate follows from summing over all divisors of $n,m$ and using the fact that $(s_n,r_m)\leq (n,m)$.

For the second part $\sum_{k=M+1}^{\infty}$, we use the fact that $|\sin(x)|\leq 1$ and obtain an absolute constant $C'>0$ with the following property
\[
\sum_{k=M}^{\infty}\frac{|\sin(2\pi\epsilon_n k)c_n(k)\sin(2\pi\epsilon_m k)c_m(k)|}{k^2}\leq  \frac{C'}{M}d(n)d(m)(n,m).
\]
We can now set $M= d(n)d(m)(n,m)n^4m^4$. We assume that $\epsilon_n\epsilon_m\neq 0$ otherwise $\lambda(A_n\cap A_m)=0$ and there is nothing to show. Then we see that $\log M\leq 10\log n+10\log m$. In particular, if $n,m\leq N$ then $\log M\leq 20\log N$. We also see that
\[
\frac{d(n)d(m)(n,m)}{M}=\frac{1}{n^4m^4}.
\]
Then, there exists an absolute constant $C''>0$ such that the following holds:
\begin{eqnarray*}
\lambda(A_n\cap A_m)\leq  4\epsilon_n\epsilon_m\phi(n)\phi(m)&+&C'' \min\{\epsilon_n,\epsilon_m\}d(n)d(m)(n,m)(10\log n\\
&+&10\log m)+C'\frac{1}{n^4m^4}.
\end{eqnarray*}
We can now use theorem \ref{T1} and lemma \ref{L3} to conclude the proof. First, observe that by Lemma \ref{L3} there exists a constant $C'''>0$ such that
\begin{eqnarray*}
& &\sum_{n=1}^N\sum_{m\leq n} \min\{\epsilon_n,\epsilon_m\}d(n)d(m)(n,m)(10\log n+10\log m)\\
&\leq& \sum_{n=1}^N\sum_{m\leq n} 20\epsilon_n d(n)d(m)(n,m)\log n\\
&\leq& C'''\sum_{n=1}^N \epsilon_n n d^3(n)\log^2 n.
\end{eqnarray*}
Similarly, the result holds for the sum $\sum_{m=1}^N\sum_{n<m}$ as well, therefore for a constant $C''''>0$ we have the following inequality
\begin{eqnarray*}(*)
\sum_{n,m=1}^N\lambda(A_n\cap A_m)&=&(\sum_{n=1}^{N}\sum_{m\leq n}+\sum_{m=1}^{N}\sum_{n< m})\lambda(A_n\cap A_m)\\&\leq& (\sum_{n=1}^N 2\epsilon_n\phi(n))^2+C''''\sum_{n=1}^{N} \epsilon_n n d^3(n)\log^2 n+100C'\zeta^2(4).
\end{eqnarray*}
From here Theorem \ref{ThmBy} follows. In fact, by the Borel-Cantelli lemma (theorem \ref{T1}), we see that
\begin{eqnarray*}
\lambda(\limsup_{n\to\infty} A_n)&\geq& \limsup_{N\to\infty} \frac{(\sum_{n=1}^N \lambda(A_n))^2}{\sum_{n,m=1}^{N}\lambda(A_n\cap A_m)}\\
&\geq& \limsup_{N\to\infty}\frac{ 1}{1+\frac{C''''\sum_{n=1}^{N} \epsilon_n n d^3(n)\log^2 n+100C'\zeta^2(4)}{(\sum_{n=1}^N 2\epsilon_n\phi(n))^2}}.
\end{eqnarray*}
The rightmost side of the above inequality is equal to $1$ under the condition of theorem \ref{ThmBy}.
Next, it is easy to see the following result for an absolute constant $C''''$ and for all integers $n$:
\[
\frac{n}{\phi(n)}d^3(n)\log^2 n\leq C'''' \log^2 n\exp(3\log2 \log n/\log\log n)\log\log n.
\]
We have used here the following result relating to the divisor function:
\[
\limsup_{n\to\infty} \frac{\log d(n)}{\log n/\log\log n}=\log 2.
\]
From here the proof of Corollary \ref{Thm1} concludes.

\section{Expected number of solutions}\label{Se3}
Here we refine the result of the previous section. The content of this section will be used in the final proof of Theorem \ref{ThmDimMain}. Previously, we have required that $f(n)\in [0,1/2]$ for all integers $n$. In this section we shall allow $f(n)$ to take any value in $[0,n/2)$. Care is needed regarding the interpretation when $f(n)>1/2$. The first thing to observe is that the following intervals for different $m$ such that $ (m,n)=1$ may overlap
\[
\left\{x: \left|x-\frac{m+\theta(n)}{n}\right|<\frac{f(n)}{n}\right\}.
\]
The second thing to observe is that it is now possible that
\[
\left\{x: \left|x-\frac{m+\theta(n)}{n}\right|<\frac{f(n)}{n}\right\}\cap (1,\infty)\neq\emptyset.
\]
To overcome these problems we need to consider $[0,1)$ as $\mathbb{R}/\mathbb{Z}$. For $x\in\mathbb{R}$ we use $\|x\|$ to be the following quantity
\[
\inf_{n\in\mathbb{Z}}|x+n|.
\] 
Given an approximation function $f$ such that for each integer $n\geq 2$, $f(n)\in [0,n/2)$ and inhomogeneous shift $\theta$ taking values in $[0,1/2)$. We want to study the following quantity for Lebesgue typical $x\in\mathbb{R}$,
\[
S(f,\theta,x,N)=\#\left|\left\{n,m\leq N, (n,m)=1:\left \|x-\frac{m+\theta(n)}{n}\right\|<\frac{f(n)}{n}\right\}\right|.
\]
We will prove here the following result.

\begin{thm}\label{ThmNm}
	For any $f:\mathbb{N}\to [0,\infty)$, $\theta:\mathbb{N}\to [0,1/2]$ and a positive number $\rho\in (0.5,1]$. If
	\[
	\limsup_{N\to\infty} \frac{\sum_{n=1}^N f(n)n^{-1} \phi(n)}{\exp(\log N\log\log\log N/\log\log N)}= \infty,
	\]
	then for Lebesgue almost all $x\in [0,1]$, there exist infinitely many integers $N_i(x)$ such that
	\[
	\left|S(f,\theta,x,N_i(x))-\sum_{n=1}^{N_i(x)} 2\frac{f(n)}{n}\phi(n)\right|\leq \left(\sum_{n=1}^{N_i(x)} 2\frac{f(n)}{n}\phi(n)\right)^\rho.
	\]
\end{thm}
\begin{proof}
As in section \ref{Se1} we construct the function \[g_n(x)=\sum_{1\leq a\leq n:(a,n)=1,\left \|x-\frac{a+\theta(n)}{n}\right\|<\frac{f(n)}{n} }1\] and see that
\[
S(f,\theta,x,N)=\sum_{n=1}^{N}g_n(x).
\]
We note here that $g_n(x)$ can take integer values other than $0$ and $1$. It is easy to see that
\[
\int_{0}^1 S(f,\theta,x,N) dx=\sum_{n=1}^N 2\epsilon_n\phi(n)=E_N.
\]
Now, we estimate the variance
\[
\int_0^1 |S(f,\theta,x,N)-E_N|^2 dx=\int_0^1 \sum_{n,m=1}^N g_n(x)g_m(x)dx-(E_N)^2.
\]
We need to consider the following integral
\[
\int_0^1 \sum_{n,m=1}^N g_n(x)g_m(x)dx=\sum_{n,m=1}^N \|g_ng_m\|_{L^1}.
\]
Although the functions $g_n$ are more complicated than the ones in Section \ref{Se1}, the computations of their Fourier coefficients are the same and results are unchanged. We omit the details here. Now we can use Fourier series to obtain the following equality as in the previous section
\begin{eqnarray*}
\|g_ng_m\|_{L^1}&=&4\epsilon_n\epsilon_m\phi(n)\phi(m)\\
&+&\frac{2}{\pi^2}\sum_{k=1}^{\infty}\frac{\sin(2\pi\epsilon_n k)c_n(k)\sin(2\pi\epsilon_m k)c_m(k)\cos\left(2\pi i \left(\frac{\theta(n)}{n}-\frac{\theta(m)}{m}\right) k\right)}{k^2}.
\end{eqnarray*}
The argument in the proof of theorem \ref{Thm1} allows us to see that for some constant $C'''>0$ we have 
\[
\int_0^1 |S(f,\theta,x,N)-E_N|^2 dx\leq C'''\sum_{n=1}^{N} \epsilon_n n d^3(n)\log^2 n.
\]
By the Markov inequality we see that given any sequence of positive numbers $\{\beta(n)\}_{n\in\mathbb{N}}$
\[
K(f,\theta,N,\beta):=\lambda(x:|S(f,\theta,x,N)-E_N|\geq \beta_N)\leq C'''\frac{1}{\beta^2_N}\sum_{n=1}^{N} \epsilon_n n d^3(n)\log^2 n.
\]
If  $\lim_{N\to\infty} K(f,\theta,N,\beta)=0,$ then there exist a subsequence $N_i$ such that
\[
\sum_i K(f,\theta,N_i,\beta)<\infty.
\]
For Lebesgue almost every $x$, there are only finitely many $N_i$ such that \[|S(f,\theta,x,N_i)-E_{N_i}|\geq\beta_{N_i}.\]
Now we see from the discussions in previous section that
\[
K(f,\theta,N,\beta)\leq C''' \frac{1}{\beta^2_N} E_N \log^2 N \log\log N \exp(3\log 2 \log N/\log\log N).
\]
Let us denote
\[
A_N=\frac{E_N}{ \exp(\log N\log\log\log N/\log\log N)},
\]
and suppose that $\limsup_{N\to\infty} A_N=\infty$, then we see that for $\beta_N=E^\rho_N$
\[
K(f,\theta,N,\beta)\to 0 \text{ as $N\to\infty$}.
\]
This is because of the following inequality and the fact that $2\rho-1>0$
\begin{eqnarray*}
K(f,\theta,N,\beta)&\leq& C''' \frac{1}{E^{2\rho-1}_N} \log^2 N \log\log N \exp(3\log 2 \log N/\log\log N)\\
&=& C''' \frac{1}{A^{2\rho-1}_N} \frac{\log^2 N \log\log N \exp(3\log 2 \log N/\log\log N)}{\exp((2\rho-1)\log N\log\log\log N/\log\log N)}.
\end{eqnarray*}
Hence for Lebesgue almost all $x\in \mathbb{R}$ there are infinitely many integers $N>0$ such that
\[
E_N+E^{\rho}_N \geq S(f,\theta,x,N)\geq E_N-E^{\rho}_N.
\]
In particular if $\rho<1$, then for such $x$ we see that for infinitely many coprime pairs $n,m$ the following inequality holds
\[
\left \|x-\frac{m+\theta(n)}{n}\right\|<\frac{f(n)}{n}.
\]
\end{proof}
\section{proof of theorem \ref{ThmDimMain}}
Recall Theorem (HS) in Section \ref{HDR}. We now show that $
\dim_H W(f,\theta)\geq \dim_H W_0(f,\mathbf{0}).
$ The other direction can be proved by the same argument provided in \cite{HS}, see also \cite[Lemma 1]{L}. Only for the lower bound are there some difficulties in estimating the size of the intersections $A_n\cap A_m$ by using direct number theoretic methods.

Let $f$ be any approximation function and $\theta$ be any inhomogeneous shift. As in the above theorem, for any $\alpha$, we find sets with cardinality $C_\alpha(N)$ and find the exponent $\delta(\alpha)$. Assume that $\kappa(\alpha)>0,$ otherwise there is nothing to show.

 First, we consider the case when $\delta(\alpha)>0$ and we shall show that
 \[
 \dim_H W(f,\theta)\geq \frac{1+\delta(\alpha)}{\alpha}.
 \]
 Now, for an arbitrarily small number $\sigma>0$ such that $\sigma<\delta(\alpha)$ we use the dimension function $h(x)=x^{\frac{1-\sigma+\delta(\alpha)}{\alpha}}$ in the mass transference principle (Theorem \ref{ThMT}). We see that $\epsilon_n=f(n)/n\geq 1/n^\alpha$ for a subset $C_\alpha$ of $\mathbb{N}$ such that
\[
\limsup_{N\to\infty} \frac{\#|C_\alpha\cap [1,N]|}{N^{\delta(\alpha)-0.5\sigma}}=\infty.
\]  
We see that $h(\epsilon_n)\geq \frac{1}{n^{1-\sigma+\delta(\alpha)}}$ and
\begin{eqnarray*}
& &\limsup_{N\to\infty} \frac{\sum_{n=1}^N \phi(n) h(\epsilon_n)}{\log^2 N\log\log N\exp(3\log 2 \log N/\log\log N)}\\
&\geq& \limsup_{N\to\infty}\frac{\#|C_\alpha\cap [1,N]| \frac{1}{\log\log N} \frac{1}{N^{-\sigma+\delta(\alpha)}}}{\log^2 N\log\log N\exp(3\log 2 \log N/\log\log N)}\\
&\geq& \limsup_{N\to\infty} \frac{N^{0.5\sigma}}{\log^2 N\log\log^2 N\exp(3\log 2 \log N/\log\log N)}=\infty.
\end{eqnarray*}
By Theorem \ref{ThMT}, we have
\[
\mathcal{H}^{\frac{1-\sigma+\delta(\alpha)}{\alpha}}(W(f,\theta))=\infty.
\]
This implies that for all $\sigma>0$
\[
\dim_H W(f,\theta)\geq \frac{1-\sigma+\delta(\alpha)}{\alpha}.
\]
This implies further that
\[
\dim_H W(f,\theta)\geq \frac{1+\delta(\alpha)}{\alpha}.
\]

Now we consider the case when $\delta(\alpha)=0$ and $C_\alpha(N)\to\infty$. For a positive number $\rho<1$ which can be chosen close to $1$, we consider the dimension function $h(x)=x^{\rho/\alpha}$. Assume that
\[
f(n)\neq 0\iff \epsilon_n=f(n)/n\geq \frac{1}{n^{\alpha}},
\]
and by shrinking some values of $f$ if necessary
\[
f(n)\neq 0\iff \epsilon_n=f(n)/n= \frac{1}{n^{\alpha}}.
\] 
Therefore we see that
\[
h(\epsilon_n)\neq 0\iff h(\epsilon_n) = \frac{1}{n^{\rho}}.
\]
Because $1/n^{\rho}>1/n$ we are in the situation discussed in section \ref{Se3}. Now if $\epsilon_n\neq 0$ we see that when $n$ is also large enough (see (***) in proof of Corollary \ref{Co1})
\[
\phi(n)h(\epsilon_n)\geq 0.0001\frac{n}{\log\log n} \frac{1}{n^{\rho}}\geq n^{0.5-0.5\rho}.
\]
This implies that
\[
\sum_{n=1}^N \phi(n)h(\epsilon_n)\geq N^{0.5-0.5\rho}
\]
for infinitely many $N$.
By theorem \ref{ThmNm} (with $f(n)=n^{1-\rho}$ in the statement), we see that for Lebesgue almost all $x\in \mathbb{R}$ there are infinitely many coprime pairs $n,m$ such that $f(n)\neq 0$ and
\[
\left \|x-\frac{m+\theta(n)}{n}\right\|<\frac{1}{n^{\rho}}.
\]
This is almost what we need, we want to find $x\in [0,1]$ such that there are infinitely many coprime pairs $n,m$ such that $f(n)\neq 0$ and
\[
\left |x-\frac{m+\theta(n)}{n}\right|<\frac{1}{n^{\rho}}.
\]
Let $M$ be a large integer. Consider now $x\in [M^{-1},1-M^{-1}]$. Suppose that there is a coprime pair $n,m$ such that
\[
\left \|x-\frac{m+\theta(n)}{n}\right\|<\frac{1}{n^{\rho}}.
\]
Then if $n$ is also large enough ($1/n^\rho<1/M^2$) we see that
\[
\left |x-\frac{m+\theta(n)}{n}\right|<\frac{1}{n^{\rho}}.
\]
This observation implies that for Lebesgue almost all $x\in [M^{-1},1-M^{-1}]$ there are infinitely many coprime pairs $n,m$ such that $f(n)\neq 0$ and
\[
\left |x-\frac{m+\theta(n)}{n}\right|<\frac{1}{n^{\rho}}.
\]
By letting $M\to\infty$ and using Theorem \ref{ThMT} we see that
\[
\mathcal{H}^{\rho/\alpha}(W(f,\theta))=\infty.
\]
This implies that
\[
\dim_H W(f,\theta)\geq \frac{\rho}{\alpha}.
\]
Now we can choose $\rho$ arbitrarily close to $1$ and observe
\[
\dim_H W(f,\theta)\geq \frac{1}{\alpha}.
\]
Then, combining this with the theorem by Hinokuma-Shiga we see that
\[
\dim_H W(f,\theta)\geq \dim_H W_0(f,\mathbf{0}).
\]
\section{proof of theorem \ref{ThDim}}
We now try to directly estimate the following sum
\[
\sum_{n,m=1}^{N}\lambda(A_n\cap A_m).
\]
\begin{thm}\label{Thmain}
	Let $f, \theta, \epsilon_n$ be as mentioned before. Then there is a constant $C>0$ such that for all integer $N>0$
	\[
	\sum_{n,m=1}^{N}\lambda(A_n\cap A_m)\leq C\left(\max_{n\in [1,N]}\epsilon_n^{0.5}\phi(n)\right)^2\log^5 N+(\sum_{n=1}^{N}2\epsilon_n\phi(n))^2
	\]
\end{thm}
\begin{proof}
	By the arguments in Section \ref{Se1} we see that
\[
\sum_{n,m=1}^{N}\lambda(A_n\cap A_m)\leq (\sum_{n=1}^{N}2\epsilon_n\phi(n))^2+\frac{2}{\pi^2}\sum_{k=1}^{\infty}\frac{1}{k^2}(\sum_{n=1}^{N}|\sin (2\pi\epsilon_n k)c_n(k)|)^2.
\]
Since $|\sin(x)|\leq 1$, for any $\alpha\in [0,1]$ we have
\[
|\sin(x)|\leq |\sin(x)|^\alpha\leq |x|^\alpha.
\]
The basic strategy is again to split the sum with respect to $k$, say,
\[
\sum_{k=1}^{\infty}=\sum_{k=1}^{M}+\sum_{k=M+1}^{\infty},
\]
for a later determined integer $M>0$. For convenience, we make the following notation:
\[
I=\sum_{k=1}^{M},
\]
\[
II=\sum_{k=M+1}^{\infty}.
\]
Then for part $I$ we use the estimate $|\sin(x)|\leq |\sin(x)|^{0.5}\leq |x|^{0.5}$
\begin{eqnarray*}
I&=&\frac{2}{\pi^2}\sum_{k=1}^{M}\frac{1}{k^2}\left(\sum_{n=1}^{N}|\sin (2\pi\epsilon_n k)c_n(k)|\right)^2\\
&\leq& \frac{4}{\pi}\sum_{k=1}^{M}\frac{1}{k^2}\left(\sum_{n=1}^{N}\epsilon_n^{0.5} k^{0.5} |c_n(k)|\right)^2\\
&\leq& \frac{4}{\pi}\sum_{k=1}^{M}\frac{1}{k}\left(\sum_{n=1}^{N}\epsilon_n^{0.5}\phi(n) \frac{|c_n(k)|}{\phi(n)}\right)^2\\
&\leq& \frac{4}{\pi}\sum_{k=1}^{M}\frac{\left(\max_{n\in [1,N]}\epsilon_n^{0.5}\phi(n)\right)^2}{k}\left(\sum_{n=1}^{N} \frac{|c_n(k)|}{\phi(n)}\right)^2.
\end{eqnarray*}
By lemma \ref{L1},\ref{L2}, we see that for a constant $C_1>0$
\[
I\leq C_1 \left(\max_{n\in [1,N]}\epsilon_n^{0.5}\phi(n)\right)^2 \log^2 N \log^3 M,
\]
where $\log^2 N$ comes from lemma \ref{L1} and $\log^3 M$ comes from lemma \ref{L2}. For $II$ we use the trivial bound $|\sin(x)|\leq 1$ and see that
\begin{eqnarray*}
II&=&\frac{2}{\pi^2}\sum_{k=M+1}^{\infty}\frac{1}{k^2}\left(\sum_{n=1}^{N}|\sin (2\pi\epsilon_n k)c_n(k)|\right)^2\\
&\leq& \frac{2}{\pi^2}\sum_{k=M+1}^{\infty}\frac{1}{k^2}\left(\sum_{n=1}^{N}|c_n(k)|\right)^2\\
&\leq& \frac{2}{\pi^2}\sum_{k=M+1}^{\infty}\frac{1}{k^2} N^4\leq C_2 \frac{N^4}{M}
\end{eqnarray*}
for another constant $C_2>0$. Note that in above inequalities we used the fact
\[
|c_n(k)|\leq \phi(n)\leq n.
\]
With some careful analysis we can replace the $N^4$ with $N^3$, but there is no essential difference as we shall see.
Now we choose $M=N^5$. The following estimate holds for a suitable constant $C>0$
\begin{eqnarray*}
I+II&\leq& 125C_1 \left(\max_{n\in [1,N]}\epsilon_n^{0.5}\phi(n)\right)^2 \log^2 N \log^3 N+C_2\frac{1}{N}\\
&\leq& C\left(\max_{n\in [1,N]}\epsilon_n^{0.5}\phi(n)\right)^2 \log^2 N \log^3 N.
\end{eqnarray*}
From here the result of this theorem follows.
\end{proof}

We can now prove theorem \ref{ThDim}:
\begin{proof}[Proof of Theorem \ref{ThDim} using Theorem \ref{Thmain}]

	By theorem \ref{Thmain} we see that for a constant $C>0$ such that
	\begin{eqnarray*}
	\sum_{n,m=1}^{N}\lambda(A_n\cap A_m)&\leq& C\left(\max_{n\in [1,N]}\epsilon_n^{0.5}\phi(n)\right)^2\log^5 N+(\sum_{n=1}^{N}2\epsilon_n\phi(n))^2,
	\end{eqnarray*}
    we have
    \begin{eqnarray*}
    \frac{(\sum_{n=1}^{N}2\epsilon_n\phi(n))^2}{\sum_{n,m=1}^{N}\lambda(A_n\cap A_m)}&\geq& \frac{(\sum_{n=1}^{N}2\epsilon_n\phi(n))^2}{C\left(\max_{n\in [1,N]}\epsilon_n^{0.5}\phi(n)\right)^2\log^5 N+(\sum_{n=1}^{N}2\epsilon_n\phi(n))^2}\\
    &\geq& \frac{1}{C \left(\max_{n\in [1,N]}\epsilon_n^{0.5}\phi(n)\right)^2\frac{\log ^{5} N}{(\sum_{n=1}^{N}2\epsilon_n\phi(n))^2}+1}.
    \end{eqnarray*}
We can then apply the following condition for $h(x)=x$
\[
\limsup_{N\to\infty} \frac{\sum_{n=1}^N \phi(n)h(\epsilon_n)}{\log^{2.5} N\left(\max_{n\in [1,N]} h(\epsilon_n)^{1/2}n\right)}=\infty,
\]
and obtain
\begin{eqnarray*}
 \limsup_{N\to\infty}\frac{(\sum_{n=1}^{N}2\epsilon_n\phi(n))^2}{\sum_{n,m=1}^{N}\lambda(A_n\cap A_m)}&\geq&1. 
\end{eqnarray*}
The conclusion of this theorem holds for the special dimension function $h(x)=x$. For general dimension functions, we can combine the special case and the mass transference principle(Theorem \ref{ThMT}) to concludes the proof.
\end{proof}

\section{Further discussions}\label{sefurther}
\subsection{Rigidity of the Hausdorff dimension}
Our result Theorem \ref{ThmDimMain} shows that the Hausdorff dimensions of sets of well approximation numbers stay unchanged under inhomogeneous shifts and dropping non-reduced fractions. We guess that this phenomena should hold in general. In order to formulate the problem we consider the following general Diophantine approximation system.

\begin{defn}
	Given any integer $n$, let $B_n$ be a subset of $\{0,1,\dots, n-1\}$. For any approximation function $f$ and inhomogeneous shift $\theta$, define
	\[
	W_B(f,\theta)=\left\{x\in [0,1]: \left|x-\frac{m+\theta(n)}{n}\right|\leq \frac{f(n)}{n} \text{ for i.m. pairs } n,m \text{ such that } m\in B_n  \right\}.
	\]
\end{defn}
Thus $W_0(.,.)$ is equal to $W_B(.,.)$ with $B_n=\{0,\dots,n-1\}$ for all integers $n$. We formulate the following two conjectures.
\begin{conj}
	For any approximation function $f$ and inhomogeneous shift $\theta$, we have the following equality
	\[
	\dim_H W_B(f,\theta)=\dim_H W_B(f,\mathbf{0}).
	\]
\end{conj}

\begin{conj}
	For any approximation function $f$ and inhomogeneous shift $\theta$, we have the following chain of inequalities
	\[
	\liminf_{n\to\infty}\frac{\log |B_n|}{\log n}\dim_H W_0(f,\mathbf{0})\leq \dim_H W_B(f,\theta)\leq \limsup_{n\to\infty}\frac{\log |B_n|}{\log n}\dim_H W_0(f,\mathbf{0}).
	\]
	In particular, if
	\[
	\lim_{n\to\infty}\frac{\log |B_n|}{\log n}=1,
	\]
	then
	\[
	\dim_H W_0(f,\mathbf{0})= \dim_H W_B(f,\theta).
	\]
\end{conj}

	\subsection{Inhomogeneous Duffin-Schaeffer problems}
The main motivation of this paper is to show that inhomogeneous metric Diophantine approximation is not too different than the homogeneous case. In fact it is a folklore conjecture that in order to prove the Duffin-Schaeffer conjecture, the homogeneous case is perhaps the hardest case. For example, in \cite{RA} it was asked whether for any inhomogeneous shift $\theta$ the following statement holds
\[
W(f,\textbf{0}) \text{ has full Lebesgue measure}\implies W(f,\theta) \text{ has full Lebesgue measure }.
\] 

There are several developments of Duffin-Schaeffer theorem in the homogeneous case. We are curious to see whether all known results about homogeneous Duffin-Schaeffer problem hold for the inhomogeneous situation as well. In particular we list below two such questions.

\begin{ques}{(See also \cite[Conjecture 1.7]{SC})}\label{q1}
	For any approximation function $f$ and inhomogeneous shift $\theta$, if the following additional condition is satisfied
	\[
	\limsup_{N\to\infty} \frac{\sum_{n=1}^{N}\frac{f(n)\phi(n)}{n}}{\sum_{n=1}^{N}f(n)}>0,
	\]
	is the following statement true
	\[
	\sum_{n=1}^{\infty} \frac{f(n)}{n}\phi(n)=\infty\implies
	W(f,\theta) \text{ has positive Lebesgue measure?}
	\]
\end{ques}

\begin{ques}
	For any approximation function $f$ that satisfies the following divergence condition
	\[
	\sum _{n=1}^\infty \frac{\phi (n) f (n)/n}{ \exp (c(\log \log n)(\log \log \log n))}= \infty,
	\]
	does the set $W(f,\theta)$ has positive Lebesgue measure for all inhomogeneous shift $\theta$?
\end{ques}

	\subsection{Cancellation of trigonometric functions}\label{Chow}

So far we have completely ignored the effect of inhomogeneous shift. In fact in \cite{RA} some dynamical shift was considered. This shed some lights on another important feature of Fourier analysis, the cancellation. Although rather technical, carefully analysis of the cancellation of trigonometric sums often provides nice results. We are curious to see whether in this case we can perform any cancellation in the main formula:
\[
\sum_{n,m=1}^N\frac{2}{\pi^2}\sum_{k=1}^{\infty}\frac{\sin(2\pi\epsilon_n k)c_n(k)\sin(2\pi\epsilon_m k)c_m(k)\cos\left(2\pi i \left(\frac{\theta(n)}{n}-\frac{\theta(m)}{m}\right) k\right)}{k^2},
\]
if in the above expression we replace the Ramanujan sums $c_n(k)$ with the full trigonometric sum
\[
\Delta_n(k)=\sum_{a\in\{0,\dots,n-1\}}e^{2\pi i ka/n}=n1_{n|k}.
\]
The last notation indicates the function equal to $n$ when $k$ is a multiple of $n$ and $0$ otherwise. In this case, in \cite[page 217, inequality (5)]{LeV}, LeVeque showed by using Fourier series method
\[
\sum_{k=1}^{\infty}\frac{\sin(2\pi\epsilon_n k)\Delta_n(k)\sin(2\pi\epsilon_m k)\Delta_m(k)\cos\left(2\pi i \left(\frac{\theta(n)}{n}-\frac{\theta(m)}{m}\right) k\right)}{k^2}\leq
2(n,m)\min\{\epsilon_n,\epsilon_m\}.
\]
Compare with our method in Section \ref{sepp}, the most significant point is that the above bound does not have any logarithmic factor. In fact, by the above inequality and the fact that for all integer $n$,
\[
\sum_{1\leq m\leq n} (n,m)=\sum_{r:r|n} r\sum_{m:m\in [1,n], (n,m)=r} 1=\sum_{r:r|n} r\phi(n/r)=n\sum_{r|n}\frac{\phi(r)}{r}\leq nd(n).
\]
With the same method as in Section \ref{sepp} we can show the following result.
\begin{thm}\label{Thmpo2}
	For any approximation function $f$ and inhomogeneous shift $\theta$, we have the following result
	\[
	\limsup_{N\to\infty} \frac{\sum_{n=1}^N f(n)}{\sqrt{\sum_{n=1}^N f(n)d(n)}}=\infty\implies
	W_0(f,\theta) \text{ has full Lebesgue measure.}
	\]
\end{thm}
Because our method in Section \ref{sepp} completely ignored the cancellation of trigonometric sums we think that by carefully performing the cancellation one can actually get rid of the logarithmic factor,
\begin{eqnarray*}
& &\sum_{n,m=1}^N\frac{2}{\pi^2}\sum_{k=1}^{\infty}\frac{\sin(2\pi\epsilon_n k)c_n(k)\sin(2\pi\epsilon_m k)c_m(k)\cos\left(2\pi  \left(\frac{\theta(n)}{n}-\frac{\theta(m)}{m}\right) k\right)}{k^2}\\
&\leq^?& C(n,m)\min\{\epsilon_n,\epsilon_m\}.
\end{eqnarray*}
Where $C>0$ is a constant and  $\leq^?$ indicates our uncertainty. If the above would be true then one could obtain the following result which is a better version of Theorem \ref{ThmBy} and a weaker version of the content of Question \ref{q1}.

\begin{conj*}
	For any approximation function $f$ and inhomogeneous shift $\theta$, we have the following result
	\[
	\limsup_{N\to\infty} \frac{\sum_{n=1}^N \phi(n)f(n)/n}{\sqrt{\sum_{n=1}^N f(n)d(n)}}=\infty\implies
	W(f,\theta) \text{ has full Lebesgue measure.}
	\]
\end{conj*}

The above argument can help us derive some new results as well. In fact, the main task is to find a good estimate for $\lambda(A_n\cap A_m)$ (see Section \ref{Se1}). We are now going to show a much weaker version of the above conjecture.
\begin{thm}\label{Thmpo}
	For any approximation function $f$ and inhomogeneous shift $\theta$, we have the following result
	\begin{eqnarray*}
	& &\sum_{n}{f(n)}=\infty \text{ and }
	\limsup_{N\to\infty} \frac{\sum_{n=1}^N \phi(n)f(n)/n}{\sum_{n=1}^N f(n)d(n)}>0\\
	&\implies&
	W(f,\theta) \text{ has positive Lebesgue measure.}
	\end{eqnarray*}
\end{thm}
\begin{proof}
We introduce the sets $\tilde{A}_t$ for integers $t$,
\[
\tilde{A}_t=\left\{x\in[0,1]| \exists 1\leq m\leq t, \left|x-\frac{m+\theta(t)}{t}\right|<\frac{f(t)}{t}\right\}.
\]
It is easy to see that $A_t\subset\tilde{A}_t$ for all integers $t$ and therefore we have that
\[
\lambda(A_n\cap A_m)\leq \lambda(\tilde{A}_n\cap \tilde{A}_m).
\]
With the help of \cite[page 217, inequality (5)]{LeV} we see that
\[
\lambda(\tilde{A}_n\cap \tilde{A}_m)\leq 4f(m)f(n)+2(n,m)\min\{\epsilon_n,\epsilon_m\}.
\]
Then, for all integer $N>0$ and a constant $C>0$, we have
\[
\sum_{n,m=1}^N \lambda(A_n\cap A_m)\leq (2\sum_{n=1}^N f(n))^2+C\sum_{n=1}^N f(n) \tilde{d}(n),
\]
where $\tilde{d}(n)$ is defined by
\[
\tilde{d}(n)=\sum_{s|n} \frac{\phi(s)}{s}.
\]
It is easy to see that $\tilde{d}(n)\leq d(n)$ and this is  enough to prove this theorem. We see that
\begin{eqnarray*}
\lambda(\limsup_{n\to\infty} A_n)&\geq& \limsup_{N\to\infty} \frac{(\sum_{n=1}^N \lambda(A_n))^2}{\sum_{n,m=1}^{N}\lambda(A_n\cap A_m)}\\
&\geq& \limsup_{N\to\infty} \frac{1}{\left(\frac{\sum_{n=1}^N f(n)}{\sum_{n=1}^N f(n)\phi(n)/n}\right)^2+C\frac{\sum_{n=1}^N f(n)d(n)}{(\sum_{n=1}^N f(n)\phi(n)/n)^2}}.
\end{eqnarray*}
In order to obtain a positive measure of the $\limsup$ set it is enough to find infinitely many integers $N_i$ and a positive number $c>0$ such that
\[
\frac{\sum_{n=1}^{N_i} f(n)\phi(n)/n}{\sum_{n=1}^{N_i} f(n)}>c,\tag{*}
\]
and
\[
\frac{(\sum_{n=1}^{N_i} f(n)\phi(n)/n)^2}{\sum_{n=1}^{N_i} f(n)d(n)}>c.\tag{**}
\]
This is almost the Duffin-Schaeffer theorem in \cite{Du} which does not require condition $(**)$. Notice that  under the condition of this theorem, $(*)$ is trivially satisfied because $f(n)\leq f(n)d(n)$ for all integers $n$. However we see that $(**)$ is satisfied even for $c=\infty$. This concludes the proof.
\end{proof}
	\subsection{Approximation functions with nice support}

	By Theorem \ref{ThmBy} and the Hardy-Ramanujan-Tur\'{a}n-Kubilius theorem on the normal order of the logarithm of the divisor function, we see that if the approximation function $f$ is supported on a large subset of $\mathbb{N}$ on which $d(n)\leq \log^{1+\epsilon} n$, then we can provide an inhomogeneous Duffin-Schaeffer type result. For a positive number $\epsilon>0$, let $A\subset\mathbb{N}$ is such that:
	\[
	a\in A\iff d(a)\leq\log^{1+\epsilon} a.
	\]
	Note that $A$ is of natural density $1$. Then, for an approximation function supported on $A$ and any inhomogeneous shift $\theta$:
	\[
	f(n)\neq 0\implies n\in A,
	\]
	and
	\[
	\limsup_{n\to\infty} \frac{\sum_{n=1}^N \phi(n)f(n)/n}{\log ^{3+3.5\epsilon}N}=\infty\implies
	W(f,\theta) \text{ has full Lebesgue measure.}
	\]
	Or in a more convenient form:
	\[
	\sum_{n=1}^{\infty}\frac{f(n)}{\log ^{3+4\epsilon} n}=\infty.\implies
	W(f,\theta) \text{ has full Lebesgue measure.}
	\]
	Note that the power $3+4\epsilon$ here is probably not optimal.

	\subsection{An example}

	We shall now discuss more about Theorem \ref{ThDim} and Corollary \ref{Co1}. A result due to Vaaler \cite{VA} says that if $f(n)=O(1/n)$, then 
	\[
	\sum_{n=1}^{\infty} \frac{f(n)}{n}\phi(n)=\infty\implies
	W(f,\mathbf{0}) \text{ has full Lebesgue measure.}
	\] 
We can provide an approximation function $f$ that does not satisfy the Duffin-Schaeffer condition $(3)$ and the extra divergence condition $(4)$ in Section \ref{Se4} nor Vaaler's condition $f(n)=O(1/n)$. To begin with, we decompose the integer set into dyadic intervals
\[
D_k=[2^k,2^{k+1}),k=0,1,\dots
\]
For each $k$, we choose an integer $m(k)$ such that
\[
\liminf_{k\to\infty} m(k)\to\infty, \sum_{k=0}^{\infty}\frac{k^{2.4}}{m(k)!}=\infty.
\]
Then in each $n\in D_k$ we assign the value $f(n)=\log^{10} n/n$ if $n$ is a multiple of $m(k)!$. Otherwise, set $f(n)=0$. It is easy to see that for large enough $n$
\[
\frac{\phi (n) f (n)}{n \exp (c(\log \log n)(\log \log \log n))}\leq \frac{1}{n\log^2 n}.
\]
Therefore, condition (\ref{BV}) is not satisfied. Next, $f(n)$ is only non zero if $n$ is a multiple of $N_n!$ for a suitable integer $N_n$ and as $n\to\infty$, $N_n\to\infty$. Then, we see that
\[
\frac{\phi(n)}{n}\leq \prod_{r \text{ prime}, r\leq N_n}\left(1-\frac{1}{n}\right)\to 0, \text{ as } n\to\infty.
\]
Hence the Duffin-Schaeffer condition (\ref{Con1}) is not satisfied. For large enough $k$ there are more than $0.5 |D_k|/m(k)!$ numbers in $D_k$ which are multiples of $m(k)!$, so we see that
\[
\sum_{n\in D_k, m(k)| n} \frac{\log^{10} n}{n \log^{7.6} n}\geq 0.5 \frac{2^k}{m(k)!} \frac{k^{10}}{2^{k+1} (k+1)^{7.6}}\geq \frac{1}{2^{10}}\frac{k^{2.4}}{m(k)!}.
\]
Here we used the fact that $k+1\leq 2k$ for all $k>1$. The conditions in Corollary \ref{Co1}. Therefore we see that $W(f,\theta)$ has full Lebesgue measure for any inhomogeneous shift $\theta.$ In particular, this holds for $\theta=\mathbf{0}.$ As we have mentioned before, this homogeneous result can be also derived from \cite[Theorem 1]{ALMTZ18}. 
\section{Acknowledgement}
HY was financially supported by the University of St Andrews. We want to thank S. Chow for providing us a draft of \cite{SC} as well as an anonymous referee for acknowledging us the research article \cite{ALMTZ18}. We also want to thank the anonymous referee(s), S. Burrell  and J. Fraser for carefully proofreading an early version of this paper.

\end{document}